\documentclass[]{article}
\usepackage[legalpaper, margin=1.2in]{geometry}

\usepackage{epstopdf,chngcntr,amsmath,amsthm,amssymb}
\usepackage[caption=false]{subfig}

\usepackage[numbers,sort&compress]{natbib}
\bibpunct[, ]{[}{]}{,}{n}{,}{,}

\theoremstyle{plain}

\theoremstyle{definition}

\theoremstyle{remark}
\newtheorem{remark}{Remark}

\newcommand{\xkh}{x_k^h}
\newcommand{\LM}{Levenberg-Marquardt }

\usepackage{algorithmic,algorithm}
\usepackage{tikz,verbatim}
\usetikzlibrary{positioning}
\usetikzlibrary{arrows,automata}
\begin{document}


\title{On the approximation of the solution of partial differential equations by artificial neural networks trained by a multilevel Levenberg-Marquardt method}

\renewcommand{\thefootnote}{\fnsymbol{footnote}}
\footnotetext[3]{INPT-IRIT, University of Toulouse and ENSEEIHT, 2 Rue Camichel, BP~7122, F-31071 Toulouse Cedex~7, France (serge.gratton@enseeiht.fr, elisa.riccietti@enseeiht.fr)}
\author{Henri Calandra,\footnote{TOTAL, Centre Scientifique et Technique Jean F\'eger, avenue de Larribau F-64000 Pau, France   (henri.calandra@total.com)} Serge Gratton\footnotemark[3], Elisa Riccietti\footnotemark[3], Xavier Vasseur\footnote{ISAE-SUPAERO, University of Toulouse, 10, avenue Edouard Belin, BP~ 54032, F-31055 Toulouse Cedex 4, France 	(xavier.vasseur@isae-supaero.fr).}}

\maketitle

\begin{abstract}
This paper is concerned with the approximation of the solution of partial differential equations by means of artificial neural networks. Here a feedforward neural network is used to approximate the solution of the partial differential equation. The learning problem is formulated as a least squares problem, choosing the residual of the partial differential equation as a loss function, whereas a multilevel Levenberg-Marquardt method is employed as a training method. This setting allows us to get further insight into the potential of multilevel methods. Indeed, when the least squares problem arises from the training of artificial neural networks, the variables subject to optimization are not related by any geometrical constraints and the standard interpolation and restriction operators cannot be employed any longer.  A heuristic, inspired by algebraic multigrid methods, is then proposed to construct the multilevel transfer operators. Numerical experiments show encouraging results related to the efficiency of the new multilevel optimization method for the training of artificial neural networks,  compared to the standard corresponding one-level procedure.  
\end{abstract}

{\bf Keywords}
Algebraic multigrid method; Artificial neural network; Multilevel optimization method; Levenberg-Marquardt method; Partial differential equation.

\section{Introduction}

The solution of partial differential equations (PDE)  modelling stationary or nonstationary physical phenomena  is a widely studied topic. Classical discretization methods employed for their solution are finite difference, finite element or finite volume methods. Typically, these methods require the solution of large linear systems, either because a fine discretization is employed (that is needed for example for the approximation of highly nonlinear solutions) or when a high dimensional problem is considered, as sampling uniformly a volume requires a number of samples which grows exponentially with its dimension. 

Multigrid methods are widely used to solve such systems, typically those arising from the discretization of linear or nonlinear elliptic partial differential equations in two or higher dimensions. They are well known for their computational efficiency and scalability \cite{book_mg,briggs,trottenberg2000multigrid}. The main idea on which multigrid methods are based is the exploitation of a hierarchy of problems, approximating the original one. Usually this hierarchy is built exploiting the fact that the underlying infinite-dimensional problem may be described at several discretization levels. The use of more levels is beneficial, as the coarse problems will generally be cheaper to solve than the original one, and it is possible to exploit the fact that such problems exhibit multiple scales of behaviour \cite{briggs}. 

The idea of making use of more levels to solve a large-scale problem has been extended also to optimization. Many multilevel optimization methods have been developed, see, e.g., \cite{fas,rmtr,mg1,mg2,mg3,mg4,mg5,mg6}. More recently a family of high-order multilevel methods have been introduced by the authors in \cite{paper_opti}.  

On the other hand, artificial neural networks (ANN) have been recently successfully employed in different fields,  such as classification and regression, image or pattern recognition, to name a few \cite{haykin2004comprehensive,krizhevsky2012imagenet,nasrabadi2007pattern}.     Their promising performance encouraged also the spread of special neural networks' hardware implementations, to decrease the computational training times and to provide a platform for efficient adaptive 
systems \cite{hardware}.  Especially, ANNs are well known for their excellent flexibility in approximating complex high-dimensional nonlinear functions.  Neural networks are then naturally suitable to approximate the solution of difficult partial differential equations. 

The use of artificial neural networks in problems involving PDEs  has attracted a lot of interest, especially in the last years. ANNs have been used for many different purposes in this field: for the numerical solution of either direct problems  \cite{ferrari2008,raissi2017c,raissi2017a,rudd2013,mishra} or of inverse problems \cite{raissi2017d,raissi2017b}, to reconstruct the equation from given data \cite{long2017,rudy2017,schaeffer2017}, or to melt with standard solution techniques such as finite element or finite difference methods, see \cite{lee1990,manevitz2004,ramuhalli2005,takeuchi1994}.

The interest in using a neural network in the solution of partial differential equations can be indeed motivated by different factors, cf.  \cite{lee1990,ritz}. First, the differential operator does not need to be discretized, and the approximate solution obtained  possesses an analytical expression. This allows not only to have an approximation to the solution at all the points of the training interval (with an accuracy usually comparable to that obtained at the training points) but also to perform extrapolation outside the interval. This approach is useful especially in case of high dimensional equations, as it allows to alleviate the curse of dimensionality, see, e.g.,  \cite{jentzen1,jentzen2,han2017} among others. More importantly, this approach provides a natural way to solve problems with nonlinear operators, with no need of linearisation. 

One of the main issues with the use of artificial neural networks is to be able to successfully train them. The training is based on the solution of an optimization problem that can be large-scale, as, like in the case of traditional techniques, in case of highly nonlinear solutions, a network with a large number of weights may be necessary to approximate the solution with a sufficient accuracy. Gradient-based methods may exhibit a slow convergence rate and it may be difficult to properly tune the learning rate. 
Recently, a new stream of research has emerged, on methods that make use of second-order derivative approximations, contributing to build the next generation of optimization methods for large-scale machine learning \cite{nocedal,bottou2018optimization}.

Inspired by all these developments, we propose the following artificial neural network based approach for the approximation of the solution of partial differential equations. We propose to use a feedforward neural network and to choose the residual of the (possibly nonlinear) differential equation as the training loss function. This gives rise to a nonlinear least squares problem. Inspired by the great efficiency of classical multigrid methods for the solution of partial differential equations, and by the recent development of second order methods in the machine learning community, we solve the training problem by a multilevel \LM method, recently introduced in \cite{paper_opti}.

The aim of our manuscript is then twofold. 

First, we want to get more insight into the potential of multilevel strategies for nonlinear optimization, investigating their application on problems that do not possess an underlying geometrical  structure. Multilevel techniques require the construction of  transfer operators. Usually standard interpolation and restriction operators, in the case of linear systems, are used also in the case of nonlinear problems. However, to employ such operators a geometrical structure of the variables is needed. In the case of a training problem, the variables subject to optimization are the weights and biases of the network, which do not possess any geometrical structure. 
To derive a multilevel method for the training problem, it is then necessary to design multilevel transfer operators differently. 

It is clear that the network possesses an important structure, even if it is not  geometric in nature. It rather presents an algebraic structure, as the one that is usually exploited in algebraic multigrid, see, e.g., \cite{ruge} and \cite[Appendix A]{trottenberg2000multigrid}. We investigate then if we could exploit this structure to define both the hierarchy of problems and the  multilevel transfer operators.  We propose a technique inspired by classical algebraic multigrid methods to do so.


Our second aim is to investigate the usefulness of multilevel methods as training methods. We inquire if the multilevel nature of the proposed solver can help to speed up the training process with respect to the one level optimization strategy. This manuscript represents a first step in this direction. Indeed, here we focus on the simplest case, that of one-layer networks. This case is interesting on his own, motivated by the Hecht-Nielsen theorem \cite{nielsen}, cf. Section \ref{sec_ANN}. From this we know that a one-layer network is capable of approximating a wide class of functions, up to any given accuracy level. The number of nodes necessary may however be really large. We experimentally show that it is not really the case for the class of problems we consider. Numerical results show indeed that with a reasonable number of nodes we can reach the desired solution accuracy. Moreover, the strategy we present to cope with the large size of the input space may further encourage the use of such networks, the extension of this to more complicated multilayer being deferred to a forthcoming paper.  


\subsection{Our contribution} We summarize here the novelties and contributions of the paper. The aim of the paper is to study the applicability of multilevel optimization methods for problems that do not possess an underlying geometrical structure, and in particular we focus on the training of artificial neural networks. Two practical questions are addressed.

\begin{itemize}
	\item \textit{How to make the approach practical and face the lack of a geometrical structure of the underlying problem?} We propose the use of a heuristic inspired by classical algebraic multigrid methods to define the hierarchy of problems and the  multilevel transfer operators. 
	\item \textit{Which is the performance of the proposed method as compared to the standard one-level version?}   We consider the approximation of a solution of a PDE by a neural network, problem lately widely addressed by the machine learning community. We show through numerical tests the gains, in term of floating points operations, arising from the use  of multilevel solvers as compared to standard one level solvers. 
\end{itemize} 

The idea of exploiting multiple scales in learning is not new, we mention for example \cite{haber2017a,stella,svm}.  Therein the multilevel structure is introduced in the model architecture, while here the training strategy is a multilevel strategy, the network's architecture being unchanged.

Moreover, to our knowledge, multilevel optimization techniques have only been applied to problems in which the hierarchy of function approximations could be built by exploiting the underlying geometrical structure of the problem at hand. Thus this work represents an improvement in the study of multilevel optimization methods.

\subsection{Organization of the paper}

The manuscript is organized as follows. In Section \ref{sec_method} we briefly review the standard \LM method, whereas in Section \ref{sec_multilevel} we describe its multilevel extension detailed in \cite{paper_opti}. Then, in Section \ref{sec_ann_pde}, we  describe the artificial neural network approximation of the solution of the partial differential equation we employ and the related least squares problem. We then discuss its solution by the multilevel solver. In particular, we introduce the heuristic we propose to build the multilevel transfer operators. Finally, in Section \ref{sec_num}, we present detailed numerical experiments related to the solution of both linear and nonlinear partial differential equations. Conclusions are drawn in Section \ref{sec_conclusion}.

\section{The \LM method}\label{sec_method}
The \LM (LM) method is an iterative procedure for the solution of least squares problems. 

Let us consider a least squares problem of the form: 
\begin{equation}\label{pb}
\min_x f(x)=\frac{1}{2}\|F(x)\|^2,
\end{equation}
with $F:\mathcal{D}\subseteq \mathbb{R}^n\rightarrow\mathbb{R}^m$, $m\geq n$, a twice continuously differentiable function.
At each iteration $k$, given the current iterate $x_k$, the objective function $f$ is approximated by the norm of an affine model of $F$, resulting in a quadratic Taylor model for $f$ with approximated Hessian matrix:
\begin{equation*}
T_k(x_k,s)=\frac{1}{2}\|F(x_k)\|^2+(J(x_k)^TF(x_k))^Ts+\frac{1}{2}s^TB_ks,
\end{equation*} 
where $J$ is the Jacobian matrix of $F$ and $B_k=J(x_k)^TJ(x_k)$ approximates the Hessian matrix $\nabla_x^2 f(x_k)$. 
This is then regularized yielding:
\begin{equation}\label{taylor_model}
m_k(x_k,s)=T_k(x_k,s)+\frac{\lambda_k}{2}\|s\|^2,
\end{equation}  
for $\lambda_k>0$ a positive value called regularization parameter.

This model is minimized (possibly approximately) to find a step $s_k$ that is used to define the new iterate $x_{k+1}=x_k+s_k$.

At each iteration it has to be decided whether to accept the step or not. This decision is based on the accordance between the decrease in the function (\textit{actual reduction}, $ared=f(x_k)-f(x_k+s_k)$) and in the model (\textit{predicted reduction}, $pred=T_k(x_k)-T_k(x_k,s_k)$):
\begin{equation}\label{rho}
\rho_k=\frac{ared}{pred}=\frac{f(x_k)-f(x_{k+1})}{T_k(x_k)-T_k(x_k,s_k)}.
\end{equation} 
If the model is a sufficiently accurate approximation to the objective function, $\rho_k$ will be close to one. Then, the step is accepted if $\rho_k$ is larger than a chosen threshold $\eta_1\in(0,1)$ and is rejected otherwise. In the first case the step is said to be \textit{successful}, otherwise the step is \textit{unsuccessful}. 

After the step acceptance, the regularization parameter is updated for the next iteration. The update is still based on the ratio \eqref{rho}. If the step is successful the parameter $\lambda$ is decreased, otherwise it is increased. 

The whole process is stopped when a minimizer of $f$ is reached. Usually, the stopping criterion is based on the norm of the gradient, i.e. given a threshold $\epsilon>0$ the iterations are stopped as soon as $\|\nabla_x f(x_k)\|<\epsilon$.

The main computational work per iteration in this kind of methods is represented by the minimization of the model \eqref{taylor_model}. This is the most expensive task, and the cost depends on the dimension of the problem. However, from the convergence theory of such methods, it is well known that it is not necessary to minimize the model exactly to obtain global convergence. 

A well-known possibility is indeed to minimize the model until the Cauchy decrease is achieved, i.e. until a fraction of the decrease provided by the Cauchy step  (the step that minimizes the model in the direction of the negative gradient) is obtained.
Here, we will consider a different kind of stopping criterion for the inner iterations, initially proposed in \cite{arc,Birgin2017}. In this case, the inner iterations (for the minimization of the model) are stopped as soon as the norm of the gradient of the regularized model becomes lower than a multiple of the squared norm of the step:
\begin{equation}\label{stopping_inner}
\|\nabla_s m_k(x_k,s_k)\|\leq \theta \|s_k\|^2,
\end{equation}
for a chosen constant $\theta>0$. 
In order to minimize the model approximately, it is possible to use a Krylov method on the system
\begin{equation*}
(B_k+\lambda_k I)s=-J(x_k)^TF(x_k), 
\end{equation*}
and stop it as soon as the inequality \eqref{stopping_inner} is satisfied. 
\vskip 5pt
The \LM procedure is sketched in Algorithm \ref{algo0}.

In the next section we will briefly review the multilevel extension of the  Levenberg-Marquardt method. This method is part of the family of methods introduced in \cite{paper_opti}, and corresponds to the case $q=1$, but with a different norm for the regularization term. In \cite{paper_opti}, if $q=1$, the regularized model is defined as 
\begin{equation}\label{model_paper_opti}
f(x_k)+\nabla_x f(x_k)^Ts+\frac{\lambda_k}{2}\|s\|^2,
\end{equation} 
where, in case of a least squares problem, $\nabla_x f(x_k)=J(x_k)^TF(x_k)$.
For a symmetric positive definite matrix $M\in\mathbb{R}^{n\times n}$ and $x\in\mathbb{R}^n$, we can define the following norm:
\begin{equation*}
\|x\|_M=x^TMx.
\end{equation*}
If we define $M_k=\displaystyle \frac{B_k}{\lambda_k}+I$, then we have $\displaystyle \frac{\lambda_k}{2}\|s\|_{M_k}^2=\frac{1}{2}s^TB_ks+\frac{\lambda_k}{2}\|s\|^2$, so that the model in \eqref{taylor_model} can be written as
$$m_k(x_k,s)=f(x_k)+\nabla_x f(x_k)^Ts+\frac{\lambda_k}{2}\|s\|_{M_k}^2,$$ corresponding to the model in \eqref{model_paper_opti}, just with a different norm for the regularization term. 

The theory presented in \cite{paper_opti} for the case $q=1$ applies for the multilevel method presented in Section \ref{sec_multilevel}, because  the $\|\cdot\|_{M_k}$ norm and the Euclidean norm are equivalent, if we assume that $\|B_k\|$ is bounded at each iteration $k$, which is a common assumption in optimization. 

\begin{algorithm}
	\caption{LM$(x_0, \lambda_0, \epsilon)$ (Standard \LM method)}
	\label{algo0}
	\begin{algorithmic}[1]
		\STATE{Given $0<\eta_1\leq\eta_2<1$,  $0<\gamma_2\leq\gamma_1< 1<\gamma_3$, $\lambda_{\min}>0$, $\theta>0$.}
		\STATE{{\bf Input:} $x_0 \in \mathbb{R}^n$, $\lambda_0 > \lambda_{\min} $, $\epsilon > 0$.}
		\STATE{$k = 0$}
		\WHILE{$\|\nabla_x f(x_k)\|>\epsilon$}
		\STATE{$\bullet$ {\bf Initialization:} Define the model $m_{k}$ as in \eqref{taylor_model}.}
		\STATE{$\bullet$ {\bf Model minimization:} Find a step $s_k$ that sufficiently reduces the model, i.e. such that \eqref{stopping_inner}  holds.}
		\STATE{$\bullet$ {\bf Acceptance of the trial point:}\label{step_acceptance0} Compute $\rho_k=\displaystyle \frac{f(x_k)-f(x_k+s_k)}{T_{k}(x_k)-T_{k}(x_k,s_k)}$.
			\IF{$\rho_k\geq \eta_1$} \STATE{$x_{k+1}=x_k+s_k$} \ELSE \STATE{$x_{k+1}=x_k$.}\ENDIF}
		\STATE{$\bullet$ {\bf Regularization parameter update:} }
		\IF{$\rho_k\geq \eta_1$} \STATE{
			$$\lambda_{k+1} =
			\bigg \{
			\begin{array}{ll}
			\max\{\lambda_{\min},\gamma_2\lambda_k\},  & \text{ if }\rho_k\geq \eta_2, \\
			\max\{\lambda_{\min},\gamma_1\lambda_k\},  & \text{ if } \rho_k< \eta_2,\\
			\end{array}
			$$}\ELSE \STATE{set 
			$\lambda_{k+1}= \gamma_3\lambda_k$.}\ENDIF
		\STATE{$k = k + 1$}
		\ENDWHILE
	\end{algorithmic}
\end{algorithm}

\section{Multilevel extension of the \LM method}\label{sec_multilevel}
We consider a least squares problem of the form \eqref{pb}.
At each iteration of the standard method, the objective function is approximated by the regularized Taylor model \eqref{taylor_model}. The minimization of \eqref{taylor_model} represents the major cost per iteration of the methods, which crucially depends on the dimension $n$ of the problem.
We want to reduce this cost by exploiting
the knowledge of alternative simplified expressions of the objective function. More specifically, we assume that we know a collection of functions $\{f^l\}_{i=0}^{l_{\max}}$ 
such that each $f^l$ is a twice-continuously differentiable function from $\mathbb{R}^{n_l}\rightarrow\mathbb{R}$  and  $f^{l_{\max}}(x)=f(x)$ for all $x\in\mathbb{R}^n$. We will also assume that, for each $i = 1,\dots,l_{\max}$, $f^l$ is â€œmore
costlyâ€ to minimize than $f^{l-1}$.

The method is recursive, so it suffices to describe the setting with two levels only. Then, for sake of simplicity, from now on we will assume that we have two approximations to our objective function $f$.  

For ease of notation, we will denote by $f^h:\mathbb{R}^{n_h}\rightarrow\mathbb{R}$ the approximation at the highest level, then $n_h=n$ and $f^h=f^{l_{\max}}$ in the notation previously used, while for $n_H<n_h$, $f^H:\mathcal{D}\subseteq\mathbb{R}^{n_H}\rightarrow\mathbb{R}$ is the approximation that is cheaper to optimize. All the quantities on the fine level will be denoted by a superscript $h$ and all the quantities on the coarse level will be denoted by a superscript $H$.
The main idea is then to use $f^H$ to construct, in the neighbourhood of the current iterate $x_k^h$ an alternative model $m_k^H$ to the Taylor model $T_k^h$ for $f^h=f$ 
\begin{equation}\label{taylor_model2}
T_k^h(x_k^h,s)=f^h(x_k^h)+\nabla_x f^h(x_k^h)^T s+\frac{1}{2}s^T B_ks.
\end{equation}
The alternative model $m_k^H$ should be cheaper to optimize than the quadratic
model $T_k^h$, and will be used, whenever suitable, to define the
step for the \LM algorithm. Of course, for $f^H$ to be useful at all in minimizing $f^h$, there should be some
relation between the variables of these two functions. We henceforth assume that there exist two full-rank linear operators $R:\mathbb{R}^{n_h}\rightarrow\mathbb{R}^{n_H}$ and $P:\mathbb{R}^{n_H}\rightarrow\mathbb{R}^{n_h}$ such that 
\begin{equation*}
\sigma P=R^T, \quad \max\{\|R\|,\|P\|\}\leq \kappa_R
\end{equation*}
for some constants $\sigma>0$ and $\kappa_R>0$. 

Let $x_{0,k}^H:=Rx_k^h$ be the starting point at coarse level. 
We define the lower level model $m_k^H$ as a modification of the coarse function $f^H$. $f^H$ is modified adding  a linear term to enforce the relation:
\begin{equation}\label{first_order_coher}
\nabla_s m_k^H(x_{0,k}^H)=R\nabla_x f^h(x_k^h).
\end{equation}
Relation \eqref{first_order_coher} crucially ensures that  the first-order behaviour of $f$ and $m_k^H$, are coherent in a
neighbourhood of $x_k^h$ and $x_{0,k}^H$. Indeed, if $s^H\in\mathbb{R}^{n_H}$ and $s^h=Ps^H$, it holds:
\begin{align*}
\nabla_x f^h(x_k^h)^Ts^h= \nabla_x f^h(x_k^h)^TPs^H=\frac{1}{\sigma} (R\nabla_x f^h(x_k^h))^Ts^H=\frac{1}{\sigma} \nabla_s m^H(x_{0,k}^H)^Ts^H.
\end{align*}
To achieve this, we define $m_k^H$ as:
\begin{equation}
m_k^H(x_{0,k}^H,s^H)=f^H(x_{0,k}^H+s^H)+(R \nabla_x f^h(x_k^h)-\nabla_x f^H(x_{0,k}^H))^Ts^H, \label{lower_model}
\end{equation}
where $\nabla_x f^h$ and $\nabla_x f^H$ are the gradients of the respective functions. At each generic iteration $k$ of our method, a step $s^h_k$ has to be computed to decrease the objective function $f$. Then, two choices are possible: the Taylor model \eqref{taylor_model2} or the lower level model  \eqref{lower_model}.
Obviously, it is not always possible to use the lower level model. For example, it may happen that $\nabla_x f^h(x_k^h)$ lies in the nullspace of $R$ and thus that $R \nabla_x f^h(x_k^h)$ is zero while $ \nabla_x f^h(x_k^h)$ 
is not. In this case, the current iterate appears to be first-order critical for $m_{k}^H$  while it is not for $f$. Using the model $m_{k}^H$ is hence potentially useful
only if $\|\nabla_s m_{k}^H(x_{0,k}^H)\| = \|R \nabla_x f(x_k^h)\|$ is large enough compared to $\|\nabla_x f(x_k^h)\|$. We therefore restrict
the use of the model $m_k^H$ to iterations where
\begin{equation}\label{go_down_condition}
\|R\nabla_x f^h(x_k^h)\|\geq \kappa_H \|\nabla_x f^h(x_k^h)\| \text{ and } \|R \nabla_x f^h(x_k^h)\|> \epsilon_H
\end{equation}
for some constant $\kappa_H\in(0,\min\{1,\|R\|\})$ and where $\epsilon\in(0,1)$  is a measure of the first-order criticality for $m_k^H$ that is judged sufficient at level $H$. Note that,
given $\nabla_x f^h(x_k^h)$ and $R$, this condition is easy to check, before even attempting to compute
a step at level $H$.

If the Taylor model is chosen, then we just compute a standard \LM step, minimizing (possibly approximately) the corresponding regularized model. If the lower level model is chosen, then we minimize its regularized counterpart 
\begin{equation}
m_{k}^H(x_{0,k}^H,s^H)+\frac{\lambda_k}{2}\|s^H\|^{2}
\end{equation} 
(possibly approximately) and get a point $x_{*,k}^H$ such that (if the minimization is successful) the model is reduced, and a step $s_k^H=x_{*,k}^H-x_{0,k}^H$. This step has to be prolongated on the fine level. Then, the step is defined as $s_k^h=Ps_k^H$.

In both cases, after the step is found, we have to decide whether to accept it or not. The step acceptance is based on the ratio:
\begin{equation}
\rho_k=\frac{f^h(x_k^h)-f^h(x_k^h+s_k^h)}{pred}.
\end{equation}
where $pred$ is defined as
\begin{itemize}
	\item $pred=m_{k}^H(Rx_k^h)-m_{k}^H(x_{*,k}^H)=m_{k}^H(Rx_k^h)-m_{k}^H(Rx_k^h,s_k^H)$  if the lower level model has been selected,
	\item $pred=T_k^h(x_k^h)-T_k^h(x_k^h,s_k^h)$ if the Taylor model has been selected. 
\end{itemize}
As in the standard form of the methods, the step is accepted if it provides a sufficient decrease in the function, i.e. if given $\eta>0$, $\rho_k\geq\eta$. 

We sketch the procedure in Algorithm \ref{algo}. As we anticipated, the procedure is recursive. For sake of clarity, allowing the possibility of having more then two levels, we denote the quantities on each level by a superscript $l$. We label our procedure MLM (Multilevel Levenberg-Marquardt). It is assumed to have at disposal a sequence of functions $\{f^l=\frac{1}{2}\|F^l\|\}_{l=1}^{l_{\max}}$ with $F^l:\mathbb{R}^{n_l}\rightarrow\mathbb{R}^{m}$ for $n_l>n_{l-1}$, with corresponding Jacobian $J^l$ and we define $B_k^l=J^l(x_k^l)^TJ^l(x_k^l)$.

\begin{algorithm}{}
	\begin{algorithmic}[1]
		\caption{MLM$(l, \mathrm{corr}^l, x_{0}^l, \lambda_{0}^l, \epsilon^l)$ (Multilevel Levenberg-Marquardt method)}\label{algo}
		\STATE{Given $0<\eta_1\leq\eta_2<1$,  $0<\gamma_2\leq\gamma_1< 1<\gamma_3$, $\lambda_{\min}>0$.}
		\STATE{{\bf Input:}  $l \in \mathbb{N}$ (index of the current level, $1 \le l \le l_{\max}$, $l_{\max}$ being the highest level), $\mathrm{corr}^l\in\mathbb{R}$ correction term to ensure coherence with upper level ($\mathrm{corr}^{l_{\max}}=0$), $x_{0}^l \in \mathbb{R}^{n_l}$, $\lambda_{0}^l > \lambda_{\min} $, $\epsilon^l > 0$.
		}
		
		\STATE{$R_l$ denotes the restriction operator from level $l$ to $l-1$, $P_l$ the prolongation operator from level $l-1$ to $l$.}
		\STATE{$k = 0$}
		\WHILE{$\|\nabla_x f^l(x_k^l)\|>\epsilon^l$}
		\STATE{$\bullet$ \label{model_choice}{\bf Model choice:}
			If $l>1$ compute $R_l\nabla_x f^l(x_k^l)$ and check \eqref{go_down_condition}. If $l=1$ or \eqref{go_down_condition} fails, go to Step \ref{step_taylor}. Otherwise, choose to go to Step \ref{step_taylor} or to Step \ref{step_lower}. }
		\STATE{$\bullet$ \label{step_taylor}  {\bf Taylor step computation:}  Define $B_k^l=J^l(x_k^l)^TJ^l(x_k^l)$ and find a step $s_k^l$ that sufficiently reduces
			\begin{equation*}
			\frac{1}{2}\|F^l(x_k^l)\|^2+(J^l(x_k^l)^TF(x_k^l))^Ts^l+\frac{1}{2}(s^l)^TB_k^ls^l+(\mathrm{corr}^l)^Ts^l+\frac{\lambda_k^l}{2}\|s^l\|^{2}.
			\end{equation*}. Go to Step \ref{step_acceptance}.}
		\STATE{$\bullet$ \label{step_lower} {\bf Recursive step computation:}
			Define 
			\begin{align*}
			&\mathrm{corr}^{l-1}=R_l \nabla_x f^l(x_k^l)- \nabla_x f^{l-1} (R_lx_k^l),\\
			&m_{k}^{l-1}(R_lx_k^{l},s^{l-1})=\frac{1}{2}\|F^{l-1}(R_lx_k^{l-1}+s^{l-1})\|^2+(\mathrm{corr}^{l-1})^Ts^{l-1}.
			\end{align*}
			Choose $\epsilon^{l-1}$ and call MLM($l-1$, $\mathrm{corr}^{l-1}$,$R_l~x_{k}^l$, $\lambda_k^l$, $\epsilon^{l-1}$)
			yielding an approximate solution $x^{l-1}_{*,k}$ of the minimization of $m_{k}^{l-1}$. Define  $s_{k}^l=P_l~(x^{l-1}_{*,k}-R_l~x_k^l)$ and $m_{k}^l(x_k^l,s^l)=m_{k}^{l-1}(R_lx_k^l,s^{l-1})$ for all $s^l=P_ls^{l-1}$.}
		\STATE{$\bullet$ \label{step_acceptance} {\bf Acceptance of the trial point:} Compute $\rho_k^{l}=\displaystyle \frac{f^l(x_k^l)-f^l(x_k^l+s_k^l)}{m_{k}^l(x_k^l)-m_{k}^l(x_k^l,s_k^l)}.$ }
		\IF{\label{step3a} $\rho_k^l\geq \eta_1$} \STATE{ $x_{k+1}^l=x_k^l+s_k^{l}$} \ELSE \STATE{\label{step3b} set $x_{k+1}^l=x_k^l$.}\ENDIF
		\STATE{$\bullet$ {\bf Regularization parameter update:} }
		\IF{$\rho_k^l\geq \eta_1$ } \STATE{ 
			$$\lambda_{k+1}^l =
			\bigg \{
			\begin{array}{ll}
			\max\{\lambda_{\min},\gamma_2\lambda_k^l\},  & \text{ if }\rho_k^l\geq \eta_2, \\
			\max\{\lambda_{\min},\gamma_1\lambda_k^l\},  & \text{ if } \rho_k^l< \eta_2\\
			\end{array}
			$$}
		\ELSE\STATE{set 
			$\lambda_{k+1}^l= \gamma_3\lambda_k^l$.}\ENDIF
		\STATE{$k = k + 1$}
		\ENDWHILE
	\end{algorithmic}
\end{algorithm}

\section{Artificial neural network based approach for the approximate solution of partial differential equations}\label{sec_ann_pde}
 In this section we describe the strategy we propose for the approximate solution of a PDE. 

Let us consider a stationary PDE written as:
\begin{subequations}\label{pde}
	\begin{align}
	D(z,u(z))&=g_1(z),\,\, z\in\Omega, \label{eq1}\\
	BC(z,u(z))&=g_2(z),\,\, z\in\partial\Omega, 
	\end{align}
\end{subequations}
where $\Omega\subset\mathbb{R}^N$, $N\geq 1$, is a connected subset, $\partial\Omega$ is the boundary of $\Omega$, $D$ is a differential operator, $BC$ is an operator defining the boundary conditions, and $g_1, g_2:\mathbb{R}^N\rightarrow\mathbb{R}$ are given functions. 
We remark that we do not need to make strong assumptions on $D$, we do not require it to be elliptic nor linear.

In the following, we describe the network's architecture we consider, how the training problem is formulated and how the MLM procedure is adapted to the solution of the specific problem. 

\subsection{ANN approximation to the PDE's solution}\label{sec_ANN}
Our approach is based on the approximation of the solution $u(z)$ of the PDE by an artificial neural network. 

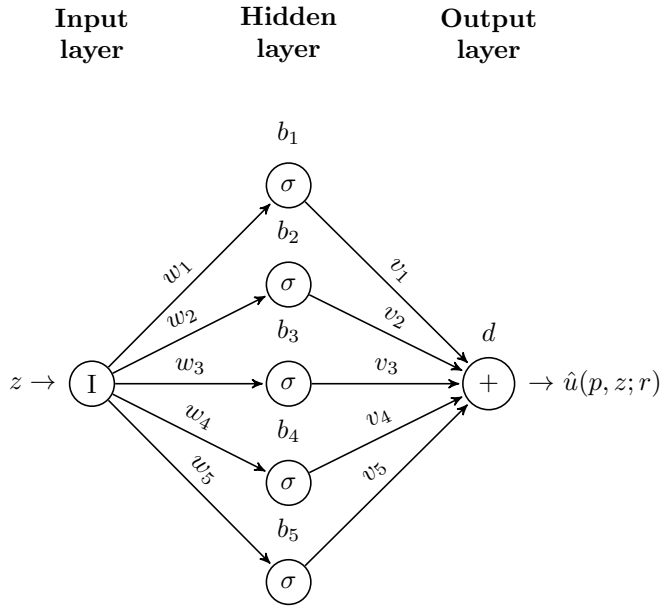
\begin{figure}
	\centering
	\begin{tikzpicture}[
	->,
	>=stealth',
	shorten >=1pt,
	auto,
	node distance=2cm,
	semithick,
	every node/.style={draw=black,circle},
	]
	\tikzstyle{annot} = [text width=4em, text centered]
	\node[label=left:{$z\rightarrow$}]         (I1)            {I};
	\node[label={$b_3$}]         (C) [right=2cm of I1]  {$\sigma$};
	\node[label={$b_4$}]         (D) [below= 0.7cm of C]  {$\sigma$};
	\node[label={$b_5$}]         (E) [below= 0.7cm of D]  {$\sigma$};
	\node[label={$b_2$}]         (B) [above= 0.7cm of C]  {$\sigma$};
	\node[label={$b_1$}]         (A) [above= 0.7cm of B]  {$\sigma$};
	
	\node[label={$d$},label=right:{$\rightarrow $ $\hat{u}(p,z;r)$}]         (O) [right=2cm of C]  {$+$};

	\node[draw=none, text width=4em, ,above of=I1, node distance=4.6cm, text centered] (hl) {\bf Input layer};
	\node[draw=none, text width=4em, ,above of=A, node distance=2cm, text centered]  {\bf Hidden layer};
	\node[draw=none, text width=4em, ,above of=O, node distance=4.6cm, text centered]  {\bf Output layer};

	\path[every node/.style={sloped,anchor=south,auto=false}]
	(I1) edge              node {$w_{2}$} (B)
	(I1) edge              node {$w_{1}$} (A) 
	(I1) edge              node {$w_{3}$} (C)          
	(I1) edge              node {$w_{4}$} (D)          
	(I1) edge              node {$w_{5}$} (E)          
	
	(B) edge              node {$v_2$\qquad} (O)
	(A) edge              node {$v_1$\qquad} (O)
	(C) edge              node {$v_3$\qquad} (O)
	(D) edge              node {$v_4$\qquad} (O)
	(E) edge              node {$v_5$\qquad} (O);
	\end{tikzpicture}
	\caption{Artificial neural network architecture with weights and biases ($r=5$, $N=1$).}
	\label{fig_weights}
\end{figure}

We assume the network to have just one hidden layer, as depicted in Figure \ref{fig_weights}. We have selected this network because the weights and biases are related in a simple way. It is therefore possible to devise a rather simple strategy to define the multilevel transfer operators.  As already discussed in the Introduction, such networks have interest on their own, motivated by the Hecht-Nielsen theorem \cite{nielsen}\footnote{For any function in $L_2[(0,1)^n]$ (i.e. square integrable on the $n$-dimensional unit cube) it exists a neural network with just one hidden layer that can approximate it, within any given accuracy.}, however the case of more hidden layers  is at the same time an interesting and a challenging problem, due to the increased nonlinearity. For this case, it is not trivial to extend the developed strategy.  We leave this as a perspective for future research. The aim of this manuscript is rather to make a first step toward a deeper understanding of the potential of multilevel techniques in  speeding up the convergence of the corresponding one level methods, when the geometry of the underlying problem cannot be exploited. 

The neural network takes the value of $z$ as input and gives an approximation $\hat{u}(p,z)$ to $u(z)$ as output, for $p\in\mathbb{R}^n$, $n=(N+2)r+1$ with $r$ the number of nodes in the hidden layer.   For the sake of simplicity, we describe our approach in the simplest case $N=1$, i.e. $z\in \mathbb{R}$. The generalization to the case $N>1$ is straightforward and is reported in Appendix A.

The network is composed of three layers in total:
one \textit{input layer} composed of just one neuron as $z\in\mathbb{R}$, that receives the value of $z$ as input; one \textit{hidden layer} with $r$ nodes, where $r$ is a constant to be fixed. A bias $b_i$, $i=1,\dots,r$, is associated with each of these nodes. All of them are connected to the input node by edges, whose weights are denoted by $w_i$, $i=1,\dots,r$, as depicted in Figure \ref{fig_weights}. The $w_i$ are called \textit{input weights}. The last layer is the \textit{output layer}, composed of just one node, as $u(z)\in \mathbb{R}$,  associated with a bias $d\in\mathbb{R}$. The output node is connected to all the nodes in the hidden layer, the corresponding weights are called \textit{output weights} and they are denoted by $v_i$, $i=1,\dots,r$.

The network is also characterized by an activation function $\sigma$, which is a given nonlinear function. Different choices are possible for $\sigma$, e.g. the sigmoid, hyperbolic tangent, logistic and softplus functions, respectively:
\begin{equation*}
\sigma(z)=\frac{e^{z}-1}{e^{z}+1},\;\;\sigma(z)=\frac{e^{2z}-1}{e^{2z}+1},\;\;\sigma(z)=\frac{e^{z}}{e^{z}+1},\;\;\sigma(z)=log(e^{z}+1).
\end{equation*}

We will denote by:
\begin{equation*}
w=[w_1,\dots,w_r]^T,\;\; v=[v_1,\dots,v_r]^T, \;\; b=[b_1,\dots,b_r]^T,\;\;p=[v,w,b,d]^T,
\end{equation*}  
the vectors of input weights, output weights and biases of the hidden nodes and the stacked vector of weights and biases, respectively.
The output of the neural network is a function of the weights and biases, but it also depends on the number of nodes $r$, which is a parameter  fixed before the training. We then denote the output as $\hat{u}(p,z;r)$, which can be expressed as 
\begin{equation}\label{exp_u}
\hat{u}(p,z;r)= \sum_{i=1}^{r}v_i \sigma(w_iz+b_i) +d, \text{ for all } z\in\Omega.
\end{equation}

The training phase consists then in the minimization of a chosen loss function, which depends on the network's output and is a function of $p$, while $z$ is constrained to the training set $\mathcal{T}$:
\begin{equation}\label{pb_ann}
\min_p \mathcal{L}(p,z)=\mathcal{F}(\hat{u}(p,z;r)),\qquad z\in\mathcal{T},
\end{equation}
for $\mathcal{F}:\mathbb{R}^{3r+1}\rightarrow\mathbb{R}$.

We choose the nonlinear residual of the equation as a  loss function, plus a penalty term with parameter $\lambda_p>0$, to impose the boundary conditions as in \cite{dissanayake1994,lagaris1998,shirvany2008}, so that for $z\in\mathcal{T}$:
\begin{equation}\label{obj_pde}
\mathcal{L}(p,z)=\frac{1}{2t}\left(\|D(z,\hat{u}(p,z;r))-g_1(z)\|^2+\lambda_p\|BC(z,\hat{u}(p,z;r))-g_2(z)\|^2\right),
\end{equation}
with $\mathcal{T}$ a training set such that $\lvert \mathcal{T}\rvert=t$, in this case a set of points $z_i$ in $\Omega$, $i=1,\dots,t$. 
\begin{remark}
	The proposed strategy does not require the discretization of operator $D$, $D(z,\hat{u}(p,z;r))$ can analytically be computed using the known expression \eqref{exp_u} of $\hat{u}(p,z;r)$.
\end{remark}

\begin{remark}
	The proposed strategy, coupled with a discretization scheme in time, can be generalized to nonstationary equations.
\end{remark}

\subsection{Solution of the training problem}\label{sec_multigrid}
Optimizing \eqref{obj_pde} may be really expensive for certain problems. If the solution $u(z)$ is highly nonlinear, a really large number of nodes $r$ may be necessary to approximate it with a sufficient accuracy. Problem \eqref{obj_pde} becomes then a large scale problem. We solve it thanks to the multilevel \LM strategy described in Section \ref{sec_multilevel}. 

We remark that, from the analysis in \cite{paper_opti}, we do not need to make strong assumptions on the operator $D$ or on the functions $g_1, g_2$ to have a convergent training method. It is just sufficient that the resulting objective function \eqref{obj_pde} has Lipschitz continuous gradient, as well as its coarse approximations. We do not need to require the operator to be elliptic nor linear.  

To be able to employ MLM in the solution of problem of the form  \eqref{pb_ann}, we need to define a strategy to build the  hierarchy of coarse problems and the multilevel transfer operators. 

Usually, when the optimization problem to be solved directly arises from the discretization of an infinite dimensional problem, the approximations to the objective function are simply chosen to be the functions arising from the discretization on coarser levels, and $P$ and $R$ are chosen to be the interpolation and the restriction operators, see \cite{briggs}. However, we remind that in our case  the variables subject to optimization are the weights and biases of the network, and not the components of the solution. Consequently, there is no geometric structure that can be exploited to construct a hierarchy, and the grids will not be real geometric grids as in classical multigrid, but rather just  sets of variables of different sizes. We have then to decide how to construct a hierarchy of sets of variables. 

Inspired by the fact that the network possesses an intrinsic algebraic structure, as the one that is typically exploited in algebraic multigrid (AMG) \cite{trottenberg2000multigrid}, we propose here a coarsening strategy based on AMG, to both explore and exploit the structure of the network, to build the hierarchy of problems and the multilevel operators. First, we give a brief overview of classical AMG techniques, and then we present the coarsening strategy we propose.

\subsubsection{Algebraic multigrid Ruge and St\"uben coarsening strategy}

AMG techniques are multilevel strategies used for the solution of linear systems of the form $Ax=b$. The goal in AMG is
to generalize the multilevel method used in geometric multigrid to target problems where the correct coarse problem is not apparent, since a geometric structure is not present. While, in geometric multigrid, a multilevel hierarchy is directly determined
from structured coarsening of the problem, in standard AMG the coarse problems, together with the transfer operators, are automatically
constructed  exploiting exclusively the information contained in the entries of the matrix. Specifically, the variables on the fine level are split into two sets, $C$ and $F$, of \textit{coarse} and \textit{fine} variables, respectively. The variables in $C$ are selected to be the variables of the coarse problem, those in $F$ are all the others. 
The coarse variables are chosen to be \textit{representative} of the fine ones, i.e. they are chosen to be the variables to which many of the remaining ones are connected. The connection among variables is only based on the entries of the matrix. There are many strategies to build this $C/F$ splitting. We decided to use the \textit{Ruge and St\"uben} strategy, one of the most classical AMG coaesening strategies \cite{brandt2000general,clees2005amg,trottenberg2000multigrid}. This is a first attempt and we do not claim that this is the best strategy to use. Other options are of course possible that could be more effective. 

The method relies on theoretical results if it is applied to a specific class of matrices \cite{trottenberg2000multigrid}, but it is commonly used also  for different problems, for which there are no theoretical guarantees.  If applied to systems arising from the discretization of simple elliptic differential operators, as the Laplace operator, this strategy is known to recover the structure exploited by geometric variants.

The splitting is built based on the notion of coupling. Two variables indexed by $i$ and $j$ are said to be \textit{coupled} if the corresponding entry of the matrix is different from zero, i.e. $a_{i,j}\neq 0$. The coupling is said to be negative if $a_{i,j}<0$, positive otherwise.  The splitting is usually made considering first negative couplings, as typically in the applications AMG is used for, the negative couplings are more than the positive ones. Then, the notion of \textit{strong negative coupling} is introduced, i.e. we say that a variable $i$ is strongly negatively coupled to another variable $j$, if  
\begin{equation}\label{strong_coupling}
-a_{i,j}\geq {\epsilon_{AMG}}\max_{a_{i,k}<0} \lvert a_{i,k} \rvert
\end{equation}
for a fixed $0<{\epsilon_{AMG}}<1$. This measure is used to construct the splitting in practice. Each $F$ variable is required to have a minimum number of its strong couplings be represented in $C$. The $C/F$ splitting is usually made choosing some first variable $i$ to become a coarse variable. Then, all variables strongly coupled to it become $F$ variables. The process is repeated until all variables have been split. In order to avoid randomly distributed $C/F$ patches, more sophisticated procedures can be designed, for example the process can be  performed in a certain order, based on a measure of importance of the variables, for more details see \cite[\S A.7.1]{trottenberg2000multigrid}. Then, also positive couplings are taken into account. After the coarsening process has been applied, a pass checks if there are strong positive $F/F$ couplings. If 
\begin{equation*}
a_{i,j}\geq {\epsilon_{AMG}} \max_{k\neq i} \lvert a_{i,k} \rvert
\end{equation*}
for some $j\neq i$, the variable $j$ is added to the set of variables strongly connected to $i$ and the variable corresponding to the largest positive coupling becomes a $C$ variable  \cite{trottenberg2000multigrid}.

Based on this splitting, the transfer operators are built to interpolate the  components of the error corresponding to the variables in $F$.  The components of the error corresponding to the variables in $C$ are transferred to the higher level by the identity operator, while the others are transformed by an interpolation formula, so that we define the $i$-th variable at fine level as: 
\begin{equation*}
x^h_i=(Px^H)_i=\begin{cases}
x_i^H\qquad\qquad\quad\text{ if }\; i\in C,\\
\sum_{k\in P_i} \delta_{i,k} x^H_k \;\;\text{ if } \;i\in F,
\end{cases}
\end{equation*}
with
\begin{equation*}
\delta_{i,k}=\begin{cases}
-\alpha_i a_{i,k}/a_{i,i} \;\text{ if }\; k\in P_i^-,\\
-\beta_i a_{i,k}/a_{i,i} \;\text{ if }\; k\in P_i^+,
\end{cases}\quad \alpha_i=\frac{\sum_{j\in N_i}a_{i,j}^-}{\sum_{k\in P_i}a_{i,k}^-},\qquad \beta_i=\frac{\sum_{j\in N_i}a_{i,j}^+}{\sum_{k\in P_i}a_{i,k}^+},
\end{equation*}
where $a_{i,j}^+=\max\{a_{i,j},0\}$, $a_{i,j}^-=\min\{a_{i,j},0\}$, $N_i$ is the set of variables connected to $i$ (i.e. all $j$ such that $a_{i,j}\neq 0$), $P_i$ the set of coarse variables strongly connected to $i$, which is partitioned in $P_i^-$ (negative couplings) and $P_i^+$ (positive couplings).
The interpolation operator, assuming to have regrouped and ordered the variables to have all those corresponding to indexes in $C$ at the beginning,  is then defined as $P=\begin{bmatrix}I\\\Delta\end{bmatrix}$ where $I$ is the identity matrix of size $\lvert C\rvert$ and $\Delta$ is the matrix such that $\Delta_{i,j}=\delta_{i,j}$, for $i=1,\dots, \lvert F\rvert$ and $j=1,\dots,\lvert C\rvert$.

\subsubsection{Algebraic coarsening strategy when solving optimization problems related to the training of artificial neural networks}

In this section, we propose a possible strategy to define both $R$ and $P$, required in the solution of problem \eqref{pb_ann}. As in classical AMG, we rely on a heuristic  strategy.
In our procedure, we have the nonlinear minimization problem \eqref{pb_ann} to solve.  At each iteration at fine level, the minimization process requires the solution of a linear system with matrix $B_k=J(x_k)^TJ(x_k)$, where $J$ is the Jacobian matrix of $F$ at $x_k$. Then, a possibility to build the $C/F$ splitting is to apply the algebraic multigrid technique we just described to $B_k$. However, we cannot apply the procedure directly to this matrix. Indeed, while in the coarsening process of standard AMG all the variables are treated in the same way, in our application the variables are coupled, cf.  \cite{clees2005amg}. We are actually optimizing with respect to triples $\{v_i, w_i, b_i\}$ of input weights, output weights and biases. The bias $d$ being a scalar is treated separately.  If no distinction is made,  a weight/bias could be removed, without the other components of the triple being removed, leading to a  network that would not be well defined.  Consequently, instead of considering the strength of connections among the variables, we will consider the strength of connections {\it{among the triples}}. We propose therefore to apply the AMG splitting to the matrix $A\in\mathbb{R}^{r\times r}$ resulting from a weighted sum of the submatrices of $B_k$ containing the derivatives of $F$ taken with respect to the same kind of variables, so that the contributions of the three different variables are not melted. More precisely, $B_k=J^T(x_k)J(x_k)$ reads: 
\begin{equation*}
B_k=\begin{bmatrix}
F_{v}^TF_v&F_v^TF_w&F_v^TF_b&F_v^TF_d\\
F_w^TF_v&F_w^TF_w&F_w^TF_b&F_w^TF_b\\
F_b^TF_v&F_b^TF_w&F_b^TF_b&F_b^TF_d\\
F_d^TF_v&F_d^TF_w&F_d^TF_b&F_d^TF_d
\end{bmatrix}, \quad F_\xi=\begin{bmatrix}
\frac{\partial F_1(x_k)}{\partial \xi_1}&\dots&\frac{\partial F_1(x_k)}{\partial \xi_r}\\
&\dots &\\
\frac{\partial F_m(x_k)}{\partial \xi_1}&\dots&\frac{\partial F_m(x_k)}{\partial \xi_r}
\end{bmatrix},\,F_d=\begin{bmatrix}
\frac{\partial F_1(x_k)}{\partial d}\\
\vdots\\
\frac{\partial F_m(x_k)}{\partial d}
\end{bmatrix},
\end{equation*}
for each variable $\xi\in\mathbb{R}^r$ and for $d\in \mathbb{R}$.

We then apply the Ruge and St\"uben splitting strategy to the following matrix:
\begin{align*}
A=\frac{F_v^TF_v}{\|F_v\|_\infty}+\frac{F_w^TF_w}{\|F_w\|_\infty}+\frac{F_b^TF_b}{\|F_b\|_\infty}.
\end{align*} 
In this way, we first obtain a $C/F$ slitting of the triples and then deduce the corresponding interpolation operator $P\in\mathbb{R}^{r\times r_c}$ with $r_c\leq r$ and the restriction operator $R\in\mathbb{R}^{r_c\times r}$ to use in the MLM method. These operators are used to project the vector $s_k^h\in\mathbb{R}^{3r}$ of fine weights and biases on the coarse level and to prolongate the coarse level step $s_k^H\in\mathbb{R}^{3r_c}$ to the fine level (in both cases omitting the contribution of the bias $d$ that is a scalar variable and is therefore left unchanged when changing levels). In both cases, the operators $P$ and $R$ are applied to each of the three components of $s_k^h, s_k^H$, corresponding to the three different kinds of variables.
In case of more than two levels, the operators $R_l$ and $P_l$, on each level, can be built with the same technique, applied to the matrix that approximates the Hessian matrix of the coarse function approximating $f$ on level $l$.

\begin{remark}
	We remark that the strategy we propose is not tailored for  \eqref{obj_pde}, it can be used for all problems of the form \eqref{pb_ann}.
\end{remark} 

\section{Numerical experiments}\label{sec_num}

In this section, we report on the practical performance of our multilevel approach. 

\subsection{Setting and parameters definition}
The whole procedure has been implemented in Julia \cite{beks:17} (version 0.6.1). We set the following values for the parameters in Algorithms \ref{algo0} and \ref{algo}, respectively: $\eta_1=0.1$, $\eta_2=0.75$, $\gamma_1=0.85$, $\gamma_2=0.5$, $\gamma_3=1.5$, $\lambda_0=0.05$ and $\lambda_{\min}=10^{-6}$. We compare the standard (one level) Levenberg-Marquardt method with a two-level Levenberg-Marquardt variant. We have decided to rely on just two levels, as in the experiments the cardinality of the coarse set of variables is much lower than that of the fine set (just a few dozens of parameters rather than more than $500$ or $1000$, depending on the problem). 

The construction of the coarse set of variables is performed through the Julia's algebraic multigrid package (AMG), that implements the classical Ruge and St\"uben method\footnote{Available at: https://github.com/JuliaLinearAlgebra/AlgebraicMultigrid.jl}.  The operators $R$ and $P$ are built just once at the beginning of the optimization procedure, using matrix $B_0$. In \eqref{go_down_condition}, we choose $\kappa_H=0.1$ and $\epsilon_H$ is equal to the tolerance chosen at the fine level. In \eqref{strong_coupling}, we set ${\epsilon_{AMG}}=0.9$, but numerical experiments highlighted that the cardinality of the coarse set does not strongly depend on this choice. For the procedure to be effective, we noticed that it is also beneficial to scale the operators $R$ and $P$, yielded by the AMG package, by their infinity norm. We choose to perform a fixed form of recursion patterns, inspired by the V-cycle of multigrid methods \cite{briggs}. Hence we alternate a fine and a coarse step, and we impose a maximum number of $10$ iterations at the coarse level.

 The linear systems arising from the minimization of the models on the fine level are solved by a truncated CGLS method \cite[\S 7.4]{bjork}, while those on the coarse level, due to the really low number of parameters, are solved  with a direct method. However, it is worth mentioning that such systems have a peculiar structure. This is similar to that of normal equations, but is fundamentally different due to the presence of the linear correction term $(R \nabla_x f^h(x_k^h)-\nabla_x f^H(x_{0,k}^H))^Ts^H$ in the right hand side, which makes it impossible to use methods for normal equations to solve the system. We refer the reader to \cite{paper_lin} for a discussion on how to exploit this structure in case of higher dimensional coarse level problems. 

\subsection{Test problems definition}
We consider both partial differential equations in one- and two-dimensions. In \eqref{eq1} we compute $g_1$ as the result of the choice of the true solution $u_T$, and we choose $u_T=u_T(z,{\nu})$, depending on a parameter ${\nu}$, that controls the oscillatory behaviour of the solution. As $\nu$ increases, the true solution becomes more oscillatory and thus harder to approximate. A network with a larger number of nodes ($r$) is thus necessary, and consequently the size of the fine problem increases. We set $\Omega=(0,1)^N$ and we impose Dirichlet boundary conditions in these experiments. We choose $\mathcal{T}$ as a Cartesian uniform grid with $h=\frac{1}{2\nu}$,  according to the Nyquist-Shannon sampling criterion \cite{shannon}, so that the cardinality of the training set is $t=2\nu+1$. Finally, in \eqref{obj_pde} we have set the penalty parameter to $\lambda_p=0.1\,t$. 

In many practical applications, as in seismology for example, a really accurate solution is not required. We stress that this setting is typical when solving partial differential equations by artificial neural networks: a high accuracy is never sought, but rather an approximate solution is looked for, see, e.g.,  \cite{mishra,raissi2017a,raissi2017b}. Therefore, we look for an approximate solution and stop the procedure as soon as $\|\nabla \mathcal{L}(p_k^h,z)\|\leq10^{-4}$ for one-dimensional problems,  $\|\nabla \mathcal{L}(p_k^h,z)\|\leq10^{-3}$ for two-dimensional problems, respectively. Indeed in the one-dimensional case, the convergence is quite fast at the beginning, allowing to quickly reach an approximate solution, while it requires far more iterations to obtain a more accurate solution for two-dimensional problems. We remark that the selected test problems are quite difficult, due to the high nonlinearity of the solutions. We outline that additional tests performed on problems with smaller value of $\nu$ (not reported here) allowed to reach a  tighter accuracy level with a smaller amount of iterations.

\begin{table}
	\centering
	\caption{List of problems considered in the following tables. We report: the equation considered, the true solution chosen and its frequency. For 2D Helmholtz equation the reference solution is computed by finite differences.}
	\begin{tabular}{c|ccc}
		\hline 
		& Equation &  $u_T(z,\nu)$&$\nu$ \\
		\hline
		Table \ref{tab_poisson}	&  $-\Delta u=g_1$ (1D) & $\cos(\nu z)$ & $20,25$ \\  
		Table \ref{tab_poisson2}	&  $-\Delta u=g_1$ (2D) &  $\cos(\nu (z_1+z_2))$ & $5,6$\\  
		Table \ref{tab:tab_helm}&	$-\Delta u-\nu^2 u=0$ (1D) & $\sin(\nu z)+\cos(\nu z)$& $5$\\
		Table \ref{tab:tab_helm2} &	$-\Delta u-\left(\frac{2\pi\nu}{c(z)}\right)^2 u=g_1$ (2D) & - & $1,2$\\
		Table \ref{tab_nonlin} (left) &$\Delta u +\sin u =g_1$ (1D)&  $ 0.1 \cos (\nu z)$& $20$ \\ 
		Table \ref{tab_nonlin} (right) &$\Delta u+e^{u}=g_1$ (2D)& $\log\left(\frac{\nu}{z_1+z_2+10}\right)$& $1$
	\end{tabular}
	\label{tab_pb}
\end{table}

\begin{table}
	\centering
	\caption{One-dimensional Poisson problem. Solution of the minimization problem (\ref{obj_pde}) with the one level Levenberg-Marquardt method (LM) and the two-level Levenberg-Marquardt method (MLM), respectively. {\tt iter} denotes the averaged number of iterations over ten simulations, ${\tt RMSE}$ the root-mean square error with respect to the true solution  and {\tt save} the ratio between the total number of floating point operations required for the matrix-vector products in LM and MLM, $r$ is the number of nodes in the hidden layer.}
	\begin{tabular}{c|ccc||ccc}
		\hline
		& & $\nu=20$ & $r=2^9$&&$\nu=25$&$r=2^{10}$\\
		\hline
		Method 	 & {\tt iter} & {\tt RMSE} &{\tt save} & {\tt iter} & {\tt RMSE} & {\tt save}  \\ 
		\hline
		LM	&  869 & $10^{-4}$ & & 1439 & $10^{-3}$& \\  
		MLM	&  507 & $10^{-4}$ &1.1-2.6-4.3 &1325 &$10^{-3}$ &1.2-1.7-2.8
	
	\end{tabular}
		\label{tab_poisson}
\end{table}

\begin{table}
	\centering
	\caption{Two-dimensional Poisson problem. Solution of the minimization problem \eqref{obj_pde} with the one level Levenberg-Marquardt method (LM) and the two-level Levenberg-Marquardt method (MLM), respectively. {\tt iter} denotes the averaged number of iterations over ten simulations, ${\tt RMSE}$ the root-mean square error with respect to the true solution and {\tt save} the ratio between the total number of floating point operations required for the matrix-vector products in LM and MLM.}
	\begin{tabular}{c|ccc||ccc}
		\hline
		& & $\nu=5$ & $r=2^{10}$&&$\nu=6$&$r=2^{11}$\\
		\hline 
		Method 	 & {\tt iter} & {\tt RMSE} &{\tt save} & {\tt iter} & {\tt RMSE} & {\tt save}  \\ 
		\hline
		LM	&  633 & $10^{-3}$& & 1213 & $10^{-3}$& \\  
		MLM	&  643 &$10^{-3}$&1.1-1.5-2.1&1016 &$10^{-3}$ &1.2-1.9-2.4 
	\end{tabular}
	\label{tab_poisson2}
\end{table}

\begin{table}
	\centering
	\caption{One-dimensional Helmholtz problem. Solution of the minimization problem \eqref{obj_pde}  with the one level Levenberg-Marquardt method (LM) and two-level Levenberg-Marquardt method (MLM), respectively. {\tt iter} denotes the averaged number of iterations over ten simulations, ${\tt RMSE}$ the root-mean square error with respect to the true solution and {\tt save} the ratio between total number of floating point operations required for the matrix-vector products in LM and MLM.}
	\begin{tabular}{c|ccc}
		\hline
		& & $\nu=5$ & $r=2^{10}$ \\
		\hline 
		Method 	 & {\tt iter} & {\tt RMSE} & {\tt save} \\ 
		\hline
		LM	&  1159 & $10^{-3}$&   \\  
		MLM	&  1250 &$10^{-3}$& 1.2-1.9-3.1
	\end{tabular}
	\label{tab:tab_helm}
\end{table}

\begin{table}
	\centering
		\caption{Two-dimensional Helmholtz problem. Solution of the minimization problem \eqref{obj_pde}  with the one level Levenberg-Marquardt method (LM) and two-level Levenberg-Marquardt method (MLM), respectively. {\tt iter} denotes the averaged number of iterations over ten simulations, ${\tt RMSE}$ the root-mean square error with respect to the  solution computed by finite differences and {\tt save} the ratio between total number of floating point operations required for the matrix-vector products in LM and MLM. With respect to the notation in Table \ref{tab_pb}, in all the tests $g_1([z_1,z_2])=(0.25< z_1<0.75 )(0.25< z_2<0.75 )$, and the velocity field $c(z)$ has been chosen as: $\bar{c}_1([z_1,z_2])=40$ (up, left); 	$\bar{c}_1([z_1,z_2])=20\,(0\leq z_1<0.5)+40\,(0.5\leq z_1\leq 1)$ (up right); $\bar{c}_2([z_1,z_2])=20\,(0\leq z_1<0.25)+40\,(0.25\leq z_2\leq 0.5)+60\,(0.5\leq z_3<0.75)+80\,(0.75\leq z_4\leq 1)$ (bottom, left); $\bar{c}_2([z_1,z_2])=0.1\sin(z_1+z_2)$ (bottom, right).}
	\begin{tabular}{c|ccc||ccc}
		\hline
		& & $\nu=1$ & $r=2^{9}$ & & $\nu=2$ & $r=2^{9}$ \\
		\hline 
		Method 	 &   {\tt iter} & {\tt RMSE}& {\tt save}&   {\tt iter} & {\tt RMSE}& {\tt save}  \\ 
		
		\hline
		LM&200&  $10^{-3}$ & &200&  $10^{-2}$ &    \\  
		MLM& 200 & $10^{-3}$ & 1.7-1.8-1.9 &200&  $10^{-2}$ & 1.7-1.8-1.9 \\
		\hline
		\hline
		& & $\nu=2$ & $r=2^{9}$ & & $\nu=2$ & $r=2^{9}$ \\
		\hline 
		Method 	 &   {\tt iter} & {\tt RMSE}& {\tt save}  &   {\tt iter} & {\tt RMSE}& {\tt save}\\ 
		
		\hline
		LM&200&  $10^{-2}$ & &200&  $5\,10^{-3}$ &    \\  
		MLM& 200 & $10^{-2}$ & 1.7-1.8-1.8 &200&  $5\,10^{-3}$ & 1.7-1.8-1.9 \end{tabular}
	\label{tab:tab_helm2}
\end{table}

The selected problems are listed in Table \ref{tab_pb}, where we report in which table the corresponding results are reported, which equation is considered, the chosen true solution and its frequency, respectively. We consider both linear and nonlinear partial differential equations and nonlinear differential equations (such as Liouville's equation). For the two-dimensional Helmholtz's equation, we have chosen different velocity fields $c(z)$: a constant function, a piecewise constant function and a sinusoidal function, cf. Table \ref{tab:tab_helm2}. For this test case, we do not choose the true solution, but impose the right-hand side. The reference solution (needed for the computation of the root mean squared error) is computed by finite differences. 

\subsection{Results of the numerical tests}
In what follows, the numerical results refer to ten simulations for different random initial guesses.  For each simulation however, the starting guess is the same for the two solvers. The first line refers to the standard (one level) Levenberg-Marquardt method (LM), while the second one to the two-level Levenberg-Marquardt method (MLM), respectively. We report the average number of iterations ({\tt iter}) over the ten simulations, the root mean squared error ({\tt RMSE}) of the computed solution evaluated on a grid of $100^N$ testing points, inside the training interval and other than the training points, with respect to the true solution evaluated at the same points and ({\tt save}), the ratio between the total number of floating point operations required for the matrix-vector products in the two methods (min-mean-max values are given), respectively. 

\begin{table}
	\centering
	\caption{Nonlinear partial differential equations. Solution of the minimization problem \eqref{obj_pde}  with the one level Levenberg-Marquardt method (LM) and two-level Levenberg-Marquardt method (MLM), respectively. {\tt iter} denotes the averaged number of iterations over ten simulations, ${\tt RMSE}$ the root-mean square error with respect to the true solution and {\tt save} the ratio between total number of floating point operations required for the matrix-vector products in LM and MLM.}
	\begin{tabular}{c|ccc||ccc}
		\hline
		& & $\nu=20$ & $r=2^{9}$& & $\nu=1$ & $r=2^{9}$\\
		\hline 
		Method 	 & {\tt iter} & {\tt RMSE} & {\tt save}& {\tt iter} & {\tt RMSE} & {\tt save}  \\ 
		\hline
		LM	&  950  & $10^{-5}$&  & 270& $10^{-3}$& \\  
		MLM	&  1444 &$10^{-5}$& 0.8-2.9-5.3 &320 &$10^{-3}$& 1.2-1.7-1.8
	\end{tabular}
	\label{tab_nonlin}
\end{table}

\vskip 5pt

As expected, we can notice that the problems become more difficult to solve as $\nu$ increases, and that the number of nodes $r$ in the neural network  has to be increased to obtain an approximate solution as accurate as for lower values of $\nu$. 

Most often, the number of iterations required by the two-level procedure is lower than that required by the one-level procedure, which is a behaviour typically observed when using classical multilevel methods. Moreover, due also to the lower dimension of the linear systems at the coarse level, the number of floating point operations required for the matrix-vector products is always {\it{considerably}} lower, even when MLM performs more iterations. The computational gain in terms of flops is on average a factor around $2$ on all the experiments, and the maximum values of the ratio is much higher for certain problems. More importantly, the quality of the solution approximation is not affected.  This is a rather satisfactory result.

We remark also that 2D tests are more difficult than 1D ones, especially those corresponding to Helmholtz's equation, which is a particularly difficult problem with Dirichlet conditions. In this case, a rough approximation is obtained in really few iterations, a maximum number of $200$ iterations is imposed, as iterating further is not useful to improve the approximation. This is already a satisfactory result. To improve the solution accuracy, a network with a more complex topology shall be more efficient, as the equation itself possess a more complex structure.

\subsection{Comments and remarks}
 These preliminary numerical results and the perspective for improvements  encourage us to investigate further on multilevel training methods.


 We finally remark that the approach introduced here is not specific to the solution of partial differential equations. It can generally be used to solve least squares problems of the form \eqref{pb_ann}, in which the solution can be expressed by a neural network. 
 
 An example that is strictly related to 1D PDEs is the solution of ODEs. These arise for example as necessary conditions of control problems. In such cases it may be interesting to find the value of the first derivative of the solution at a given point, which is straightforward with the proposed approach, as we find an analytical expression of the solution.

\section{Conclusions}\label{sec_conclusion}
We have investigated the potential of the multilevel \LM method in the solution of training problems arising from the approximation of the solution of partial differential equations by a neural network. This is chosen as representative of a class of problems in which  the variables subject to optimization are not related by any an explicit geometrical structure. To the best of our knowledge, this is also the first attempt at using multilevel optimization for the training of artificial neural networks.

We have proposed a possible heuristic strategy based on standard algebraic multigrid methods, to both explore and exploit the structure of the neural network to build a hierarchy of problems and the multilevel transfer operators. 

 This strategy has been designed for networks with one hidden layer. Considering a multilayer network is a natural and challenging extension due to the increased nonlinearity. The performance of the multilevel optimization method has been tested and compared to that of the standard one-level version. The numerical results are quite satisfactory, showing the potential of the multilevel strategy.  
We currently consider the extension of the procedure to multilayer networks as a significant perspective, as we believe that this could lead to a competitive learning method,  especially if coupled with a strategy to make the approach purely matrix-free.

\section*{Funding}
This work was funded by TOTAL.

\bibliographystyle{plain}
\bibliography{biblio_arxiv}

\appendix

\section{Appendix}
\counterwithin{figure}{section}

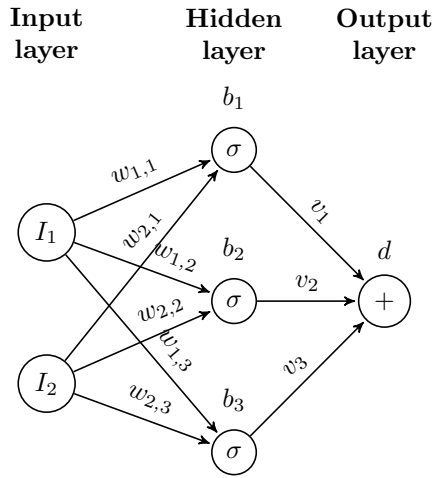
\begin{figure}
	\centering
	\begin{tikzpicture}[
	->,
	>=stealth',
	shorten >=1pt,
	auto,
	node distance=2cm,
	semithick,
	every node/.style={draw=black,circle},
	]
	\node[]         (I1)            {$I_1$};
	\node[]         (C) [below  of=I1]  {$I_2$};
	\node[label={$b_1$}]         (B) [above right=0.6cm and 2cm of I1]  {$\sigma$};
	\node[label={$b_2$}]         (A) [below of=B]  {$\sigma$};
	\node[label={$b_3$}]         (D) [below of=A]  {$\sigma$};
	\node[label={$d$}]         (O) [right of=A]  {$+$};
	
	\path[every node/.style={sloped,anchor=south,auto=false}]
	(I1) edge              node {$w_{1,1}$} (B)
	(I1) edge              node {\qquad$w_{1,2}$} (A) 
	(I1) edge              node {\qquad$w_{1,3}$} (D)          
	(C) edge              node {\qquad\,$w_{2,2}$} (A)
	(C) edge              node {$w_{2,3}$} (D)
	(C) edge              node {\qquad$w_{2,1}$} (B)
	(B) edge              node {$v_1$\qquad} (O)
	(A) edge              node {$v_2$\qquad} (O)
	(D) edge              node {$v_3$\qquad} (O);
	
	\node[draw=none, text width=4em, ,above of=I1, node distance=2.6cm, text centered] (hl) {\bf Input layer};
	\node[draw=none, text width=4em, ,above of=A, node distance=3.5cm, text centered]  {\bf Hidden layer};
	\node[draw=none, text width=4em, ,above of=O, node distance=3.5cm, text centered]  {\bf Output layer};
	\end{tikzpicture}
	\caption{Artificial neural network architecture with weights and biases (Case of $r=3$ and $N=2$).}
	\label{fig:fig_weights2D_appendix}
\end{figure}

In Section \ref{sec_ANN}, we have considered the simplest case of network's architecture, corresponding to $N=1$. However, our method is also applicable when $N>1$, as shown in Section \ref{sec_num}. Here, we describe this extension (still considering the case of just one hidden layer). The network is still composed of three layers in total, but  the input layer is composed of $N$ nodes, one for each component of the input. The input nodes are connected to all the nodes in the hidden layer. Instead of having just $r$ input weights, we have $N$ groups of $r$ input weights, $\{w_{i,1},\dots,w_{i,r}\}$, $i=1,\dots,N$, where the weights in group $i$ are located on the edges departing from the $i$-th input node.
The right part of the network remains the same, with one output node and $r$ output weights. 
The input space of the objective function $\mathcal{F}$ in \eqref{pb_ann} will be $\mathbb{R}^{(N+2)r+1}$ rather than $\mathbb{R}^{3r+1}$. 
As an example, we depict the case $N=2$ in Figure \ref{fig:fig_weights2D_appendix}. 

The training procedure is exactly the same, with the only difference that \linebreak 
$p=[v, w_1,\dots, w_N, b, d]\in \mathbb{R}^{(N+2)r+1}$, with $w_i=[w_{i,1},\dots,w_{i,r}]^T$, $i=1,\dots, N$. In the construction of the matrix required for the application of the AMG coarsening strategy, the contribution of the input weights is now represented by $N$ submatrices, corresponding to the derivatives with respect to couples $(w_i,w_i)$, $i=1,\dots,N$:
\begin{align*}
A=&\frac{F_v^TF_v}{\|F_v\|_\infty}+\sum_{i=1}^{N}\frac{F_{w_i}^TF_{w_i}}{\|F_{w_i}\|_\infty}+\frac{F_b^TF_b}{\|F_b\|_\infty}.
\end{align*} 

\end{document}